\documentclass[12pt]{article}
\usepackage{graphicx}
\usepackage{graphics}
\usepackage{amsthm}
\usepackage{amssymb, amsmath}
\title{Bernoulli identities, zeta relations, determinant expressions, Mellin transforms,
and representation of the Hurwitz numbers}
\author{Mark W. Coffey\\
Department of Physics\\
Colorado School of Mines\\
Golden, CO  80401\\
USA\\
mcoffey@mines.edu\\}
\date{January 2, 2016}
\begin{document}
\maketitle
\baselineskip=25 pt
\begin{abstract}

The Riemann zeta identity at even integers of Lettington, along with his other Bernoulli and
zeta relations, are generalized.  Other corresponding recurrences and determinant relations are
illustrated.  Another consequence is the application to sums of double zeta values. 
A set of identities for the Ramanujan and generalized Ramanujan polynomials is presented.
An alternative proof of Lettington's identity is provided, together with its generalizations to the
Hurwitz and Lerch zeta functions, hence to Dirichlet $L$ series, to Eisenstein series, and to 
general Mellin transforms.

The Hurwitz numbers $\tilde{H}_n$ occur in the Laurent expansion about the origin of a certain Weierstrass $\wp$ function for a square lattice, and are highly analogous to the Bernoulli numbers.  An integral
representation of the Laurent coefficients about the origin for general $\wp$ functions, and for these numbers in particular, is presented.  
As a Corollary, the asymptotic form of the Hurwitz numbers is determined.  In addition,
a series representation of the Hurwitz numbers is given, as well as a new recurrence.

\end{abstract}
 
\medskip
\baselineskip=15pt
\centerline{\bf Key words and phrases}
\medskip 
Bernoulli number, Bernoulli polynomial, Riemann zeta function, Euler number, Euler polynomial,
alternating zeta function, double zeta values, Hurwitz zeta function, Lerch zeta function, polygamma function, Ramanujan polynomial, Bernoulli relations, zeta identities, Eisenstein series, recurrence, Hessenberg determinant, integral representation, Mellin transform, Hurwitz numbers 

\medskip
\noindent
{\bf 2010 MSC numbers}
\newline{11B68, 11C20, 11M06}  

\baselineskip=25pt

\pagebreak
\centerline{\bf Introduction and statement of results}

\medskip
Let $\zeta(s)$ denote the Riemann zeta function and $B_n(x)$ the $n$th degree Bernoulli polynomial,
such that $B_n=B_n(0)=(-1)^{n-1}n\zeta(1-n)$ is the $n$th Bernoulli number (e.g., \cite{carlitz,lehmer}).
This relation is extended, for instance when $n$ is even, to
$B_n(x)=-n\zeta(1-n,x)$ where $\zeta(s,x)=\sum_{n=0}^\infty (n+x)^{-s}$ (Re $s>1$) is the Hurwitz zeta function.  There is the well known relation (explicit evaluation) 
$$\zeta(2m)={{(-1)^{m+1}} \over {(2m)!}}2^{2m-1}\pi^{2m}B_{2m}. \eqno(1.1)$$
The Bernoulli numbers have the well known exponential generating function
$${x \over {e^x-1}}=\sum_{n \geq 0} {B_n \over {n!}}x^n, ~~~~|x|<2\pi,$$
while the ordinary generating series $\beta(x)=\sum_{n\geq 0}B_nx^n=1-x/2+x^2/6-x^4/30+x^6/42+\ldots$
is divergent.  However the continued fractions for the generating series $\sum_{n \geq 1}B_{2n}(4x)^n$,
$\beta(x)$, $\sum_{n \geq 1}(2n+1)B_{2n}x^n$, and $\sum_{n \geq 1}(4^n-2^n)|B_n|x^{2n-1}/n$ are
convergent \cite{frame}.

In \cite{mcl1} was presented a new identity for the Riemann zeta function at even integers,
$$\zeta(2j)=(-1)^{j+1}\left[{{j\pi^{2j}} \over {(2j+1)!}}+\sum_{k=1}^{j-1}{{(-1)^k\pi^{2j-2k}} \over
{(2j-2k+1)!}}\zeta(2k)\right]. \eqno(1.2)$$
This is hardly an isolated identity, and we show how to systematically derive many related ones.

As an illustration, we have the identities given in the following.
{\newline \bf Theorem 1}. (a) 
$$(3^{1-2j}-1)\zeta(2j)={{(-1)^j j (2\pi)^{2j}} \over {(2j)!3^{2j-1}}}
+(-1)^{j+1}(2\pi)^{2j}\sum_{m=0}^j {{(-1)^{m+1}} \over {(2j-2m)!2^{2m-1}3^{2j-2m}\pi^{2m}}}\zeta(2m),$$
(b)
$$(4^{1-2j}-2^{1-2j})\zeta(2j)={{(-1)^j j (2\pi)^{2j}} \over {(2j)!4^{2j-1}}}
+(-1)^{j+1}(2\pi)^{2j}\sum_{m=0}^j {{(-1)^{m+1}} \over {(2j-2m)!2^{2m-1}4^{2j-2m}\pi^{2m}}}\zeta(2m),$$
(c)
$$(6^{1-2j}-3^{1-2j}-2^{1-2j}+1)\zeta(2j)={{(-1)^j j (2\pi)^{2j}} \over {(2j)!6^{2j-1}}}
+(-1)^{j+1}(2\pi)^{2j}\sum_{m=0}^j {{(-1)^{m+1}} \over {(2j-2m)!2^{2m-1}6^{2j-2m}\pi^{2m}}}\zeta(2m),$$
(d)
$$-{j \over 12^{2j-1}}(1+7^{2j-1})+\sum_{m=0}^j {{(2j)!} \over {(2j-2m)!}}{{(-1)^{m+1}} \over {2^{2m-1}
\pi^{2m}}}(1+7^{2j-2m}){{\zeta(2m)} \over 12^{2j-2m}}$$
$$={{(6^{1-2j}-3^{1-2j}-2^{1-2j}+1)(2j)!(-1)^{j+1}} \over {2^{4j-1}\pi^{2j}}}\zeta(2j),$$
and (e)
$$-{j \over 8^{2j-1}}(1+5^{2j-1})+\sum_{m=0}^j {{(2j)!} \over {(2j-2m)!}}{{(-1)^{m+1}} \over {2^{2m-1}
\pi^{2m}}}(1+5^{2j-2m}){{\zeta(2m)} \over 8^{2j-2m}}={{(4^{1-2j}-2^{1-2j})(2j)!(-1)^{j+1}} \over
{2^{4j-1}\pi^{2j}}}\zeta(2j).$$
Moreover, in a separate section we give an alternative proof of (1.2), and its extensions to
$\zeta(s,a)$ and so to Dirichlet $L$-functions.  

Now define functions
$$\phi_j(s)={{\zeta(s)} \over j^s}, ~~~~~~\theta_3(s)=\left(1-{1 \over 3^{s-1}}\right)\zeta(s),$$
$$\theta_4(s)=\left({1 \over 2^{s-1}}-{1 \over 4^{s-1}}\right)\zeta(s), ~~
\theta_6(s)=\left(-1+{1 \over 2^{s-1}}+{1 \over 3^{s-1}}-{1 \over 6^{s-1}}\right)\zeta(s).
\eqno(1.3)$$
Theorem 1 and similar results have numerous implications for these and other functions.
In particular, Theorem 1(a) leads to Theorem 2(b).

{\bf Theorem 2}.  (a)
$$4 {j^{2s} \over 2^{2s}}\phi_j(2s)={{\pi^{2s}(2s-1)} \over {(2s+1)!}}
+\sum_{n=1}^{s-1} {{(-1)^{s-n}\pi^{2n}} \over {(2n+1)!}}4 {j^{2(s-n)} \over 2^{2(s-n)}}\phi_j(2s-2n),$$
(b) 
$$\theta_3(2j)=(-1)^j\pi^{2j}{2^{2j} \over 3^{2j}}{{(1-3j)} \over {(2j)!}}
+2\sum_{m=0}^{j-1}{{(-1)^{m-1}} \over {(2m)!}}\left({{2\pi} \over 3}\right)^{2m}\left(1-{3 \over 3^{2(j-m)}}\right)^{-1}\theta_3(2j-2m),$$
(c)
$$\theta_4(2j)=(-1)^j\pi^{2j}{1 \over 2^{2j}}{{(1-4j)} \over {(2j)!}}
+2\sum_{m=0}^{j-1}{{(-1)^{m-1}} \over {(2m)!}}\left({{\pi} \over 2}\right)^{2m}\left({2 \over 2^{2(j-m)}}-{4 \over 4^{2(j-m)}}\right)^{-1}\theta_4(2j-2m),$$
and (d)
$$\theta_6(2j)=(-1)^j\pi^{2j}{1 \over 3^{2j}}{{(1-6j)} \over {(2j)!}}
+2\sum_{m=0}^{j-1}{{(-1)^{m-1}} \over {(2m)!}}\left({{\pi} \over 3}\right)^{2m}\left(-1+{2 \over 2^{2(j-m)}}+{3 \over 3^{2(j-m)}}-{6 \over 4^{2(j-m)}}\right)^{-1}$$
$$~~~~~~~~~~~~~~~~~~~~~~~~~~~~~~~~~~~~~~~~~~~~~~~~~~~~~~~~~~~~~~~~~~~~~~~~~~~~~\times \theta_6(2j-2m).$$


We next recall the definition of a half-weighted minor corner layered determinant 
$\Psi_s(\vec{h},\vec{H})$, using the vectors $\vec{h}=(h_1,h_2,h_3,\ldots)$ and $\vec{H}=(H_1,H_2,H_3,\ldots)$:
$$\Psi_s(\vec{h},\vec{H})=(-1)^s\left|\begin{array}{cccccc}
H_1 & 1   & 0   & 0 & \ldots & 0 \\
H_2 & h_1 & 1   & 0 & \ldots & 0 \\
H_3 & h_2 & h_1 & 1 & \ldots & 0 \\
\vdots & \vdots & \vdots & \vdots &\ddots & \vdots \\
H_{s-1} & h_{s-2} & h_{s-3} & h_{s-4} & \ldots & 1\\
H_s & h_{s-1} & h_{s-2} & h_{s-3} & \ldots & h_1 \\ \end{array}\right|.$$

{\bf Theorem 3}.  (a) Define vectors $\vec{u}$ and $\vec{U}_2$ with entries $u_s=1/(2s+1)!$ and
$U_{2s}=(2s-1)/(2s+1)!$.  Then
$$4 {j^{2s} \over 2^{2s}}\phi_j(2s)=(-1)^s \pi^{2s}\Psi_s(\vec{u},\vec{U}_2).$$
(b) Define vectors $\vec{U}_3$ and $\vec{H}_3$ with entries $U_{3s}=1/(2s)!$ and
$$H_{3s}=-\left[{{1-3s} \over {(2s)!}}+\theta_3(2s)(-1)^{s-1}\left({{2\pi} \over 3}\right)^{-2s}\right].
$$  
Then
$$(-1)^s \left({{2\pi} \over 3}\right)^{2s}\Psi_s(\vec{U}_3,\vec{H}_3)=2\theta_3(2s)\left(1-{3 \over 3^{2s}}\right)^{-1}.$$

Let ${t \choose {d_1,d_2,\ldots,d_s}}$ denote the multinomial coefficient, the coefficient of the
expansion of the sum of $s$ terms to the $t$th power.  
\footnote{On pp.\ 22 and 23 of \cite{mcl2}, the following typographical errors occur.  The upper
index for the multinomial coefficient for the summations for $\eta(2s)$, $\phi(2s)$, and $\theta(2s)$
in Lemma 3.3 should be $t$.  In the display equation for the proof of Theorem 1.3, an ``$=$" should
be inserted after $x^{2s-2}$.}
\footnote{On p.\ 17 of \cite{mcl2}, $(-1)^k$ should be $(-1)^{s-k}$ in the summand on the right side
of the second display equation, and vice versa for the summand of the right side of the third display
equation. At the bottom of p.\ 8, $\Psi$ should read $\Psi_n$ (twice).}

{\bf Corollary 1}.  
$$4\phi_j(2s)={{(2\pi)^{2s}} \over j^{2s}}{1 \over {(2^{2s-1}-1)}}\sum_{t=1}^s \sum_{d_i \geq 0}{t \choose
{d_1,d_2,\ldots,d_s}} {{(-1)^{t+s}} \over {3!^{d_1}5!^{d_2}\cdots (2s+1)!^{d_s}}},$$
and 
$$4\theta_3(2s)=\left(1-{1 \over 3^{2s-1}}\right){{(2\pi)^{2s}} \over {(2^{2s-1}-1)}}\sum_{t=1}^s \sum_{d_i \geq 0}{t \choose
{d_1,d_2,\ldots,d_s}} {{(-1)^{t+s}} \over {3!^{d_1}5!^{d_2}\cdots (2s+1)!^{d_s}}}.$$
Here the sums are such that $d_1+d_2+\ldots+d_s=t$ and $d_1+2d_2+\ldots+sd_s=s$.

In fact, the minor corner layered determinants of \cite{mcl1,mcl2} are very special cases of lower
Hessenberg determinants.  Following the proof of Theorem 3, we recall a much more general
determinantal result.

Define, for integers $a\geq 2$ and $b \geq 1$, the double zeta values
$$\zeta(a,b)=\sum_{n=1}^\infty \sum_{m=1}^{n-1} {1 \over {n^a m^b}}.$$
From manipulating series, it follows that
$$\zeta(a,b)+\zeta(b,a)=\zeta(a)\zeta(b)-\zeta(a+b).$$
We note that any recurrence of Riemann zeta values (or of Bernoulli numbers) of the form
$\sum_j f(s,j)\zeta(2j)\zeta(2s-2j)$ for some function $f$ then leads to a sum formula for double zeta values.  For we have, for instance,
$$\sum_j f(s,j)[\zeta(2j,2s-2j)+\zeta(2s-2j,2j)]=\sum_j f(s,j)\zeta(2j)\zeta(2s-2j)-\zeta(2s)\sum_j
f(s,j).$$

As an illustration we have the following.
\newline{\bf Theorem 4}. For integers $s>1$,
$$\sum_{k=1}^{s-1}{1 \over 2^{2(s-k)}}\left(1-{1 \over 2^{2k}}\right)[\zeta(2k,2s-2k)+\zeta(2s-2k,2k)]
={2^{-2s-1} \over 3}(4^s+6s-1)\zeta(2s).$$

Let $E_n(x)=\sum_{k=0}^n {n \choose k}{E_k \over 2^k}\left(x-{1 \over 2}\right)^{n-k}$, with
$E_{2n+1}=0$ and $E_k=2^kE_k(1/2)$ denote the Euler polynomial.  It has a well known exponential
generating function
$${{2e^{xt}} \over {e^t+1}}=\sum_{n=0}^\infty E_n(x) {t^n \over {n!}}, ~~~~|t|<\pi.  \eqno(1.4)$$
We recall a connection with the alternating zeta function
$$\eta(s)=\sum_{n=1}^\infty {{(-1)^{n-1}} \over n^s}=(1-2^{1-s})\zeta(s).$$
We have the evaluation
$\eta(-j)=(-1)^j E_j(0)/2$ for $j\geq 0$.  The values $E_n(0)=(-1)^nE_n(1)$ are $0$ for $n\geq 2$
and $E_0(0)=1$, while otherwise (\cite{nbs}, p. 805) $E_n(0)=-2(n+1)^{-1}(2^{n+1}-1)B_{n+1}$ for
$n \geq 1$.  Many results analogous to Theorem 1 are possible, and the following provides a brief
example.

{\bf Theorem 5}.  For integers $n>0$,
$$-{1 \over {2n}}(1-3^{1-2n})(2^{2n}-1)B_{2n}=\sum_{m=0}^{n-1}{{2n-1} \choose {2m}}{E_{2m} \over {(2m)!
}}\left(-{1 \over 6}\right)^{2(n-m)-1}.$$

{\bf Theorem 6}. For integers $j\geq 0$,
$$2\sum_{m=0}^j(1-2^{1-2m})(1-2^{2(m-j)+1})\zeta(2m)\zeta(2j-2m)=-(1-2j)\zeta(2j).$$


Theorem 1.3 of \cite{mcl2} gives equivalent forms of the recurrence
$$\zeta(2s+2)={2 \over {2^{2s+2}-1}}\sum_{k=0}^{s-1} (2^{2k+2}-1)\zeta(2s-2k)\zeta(2k+2). \eqno(1.5)$$
Hence, if we introduce functions $\tilde{\theta}_j(s)=(1-j^{-s})\zeta(s)$, we may write
$${{(2^{2s+2}-1)} \over {(1-j^{-(2s+2)})}}\tilde{\theta}_j(2s+2)=2j^{2(s+1)}\sum_{k=0}^{s-1}
{{(2^{2k+2}-1)} \over {(j^{2(k+1)}-1)}}\phi_j(2s-2k)\tilde{\theta}_j(2k+2).$$
In fact, (1.5) is proved in \cite{basu} (p. 406) as Theorem 3.4.
In addition, we may note that besides the well known relation (e.g., \cite{underwood,williams})
$$\zeta(2n)={2 \over {2n+1}}\sum_{k=1}^{n-1} \zeta(2k)\zeta(2n-2k),$$
Williams \cite{williams} some time ago proved
$$\sum_{k=1}^n{\cal L}(2k-1){\cal L}(2n-2k+1)=\left(n-{1 \over 2}\right)(1-2^{-2n})\zeta(2n),$$
where ${\cal L}$ is the function
$${\cal L}(s)=\sum_{k=0}^\infty {{(-1)^k} \over {(2k+1)^s}}=2^{-s}\Phi\left(-1,s,{1 \over 2}\right)
=4^{-s}\left[\zeta\left(s,{1 \over 4}\right)-\zeta\left(s,{3 \over 4}\right)\right],$$
and $\Phi$ is the Lerch zeta function (see Theorem 9).


We point out that {\it each} part of Theorem 1 permits the identification of candidate pseudo-characteristic polynomials for the Riemann zeta and other related functions.  
For instance, if we separate the $m=0$ and $m=j$ terms of the sum of the right side of
Theorem 1(a), we obtain 
$$3(3^{-2j}-1)\zeta(2j)=(-1)^j{{(3j-1)} \over {(2j)!}}\left({{2\pi} \over 3}\right)^{2j}
+\sum_{m=1}^{j-1} {{(-1)^m 2^{2m+1}\pi^{2m} \zeta(2j-2m)} \over {(2m)!3^{2m}}}.$$
We then put
$$p_j^{(3)}(x)={1 \over {3(3^{-2j}-1)}}\sum_{m=1}^{j-1}{{(-1)^m2^{2m+1}\pi^{2m}} \over {(2m)!3^{2m}}}x^{2m}$$
and
$$z_j^{(3)}(x)=(-1)^j{{(3j-1)} \over {(2j)!3(3^{-2j}-1)}}\left({{2\pi} \over 3}\right)^{2j}
+p_j^{(3)}(x).$$
We recognize that
$$p_j^{(3)}(x)={2 \over {3(3^{-2j}-1)}}\left[-1+\cos\left({{2\pi} \over 3}x\right)-\sum_{m=j}^\infty
{{(-1)^m} \over {(2m)!}}\left({{2\pi} \over 3}x\right)^{2m}\right].$$
It then appears that inequalities of the following sort hold, where $k$ is either $2s$ or $2s-1$.
\newline{\bf Conjecture 1}.  For integers $s \geq 4$,
$$\zeta(k)-\{\zeta(k)\}^2 \leq z_s^{(3)}(\zeta(k))\leq \zeta(k)+\{\zeta(k)\}.$$
Here $\{x\}=x-[x]$ denotes the fractional part of $x$.

We provide new identities for the Ramanujan $R_{2s+1}(z)$ and generalized Ramanujan $R_{2s}(z)$
polynomials, and for other functions.  For this we require the following definitions.
Let the numbers $B_s^*$ for $s\geq 2$ be given by
$$B_s^*=-{1 \over {s+1}}\sum_{k=0}^{s-1}{{s+1} \choose k}2^{k-s}B_k,$$
which are then such that $B_{2s}^*=B_{2s}$ and $B_{2s-1}^*=(1-2^{-2s})2B_{2s}/s$.
To these values are prepended the initial values $B_0^*=1$ and $B_1^*=1/4$.  The Ramanujan
polynomials are given by
$$R_{2s+1}(z)=\sum_{k=0}^{s+1}{{B_{2k}B_{2s+2-2k}} \over {(2k)!(2s+2-2k)!}}z^{2k},$$
while the generalized Ramanujan polynomials $Q_r(z)$ \cite{mcl2} are given by
$$Q_r(z)=\sum_{k=0}^{[(r+1)/2]} {{B_{r+1-2k}^* B_{2k}^*} \over {(r+1-2k)!(2k)!}}z^{2k}.$$
The Pochhammer symbol $(a)_j=\Gamma(a+j)/\Gamma(a)$, with $\Gamma$ denoting the Gamma function,
$\psi(z)=\Gamma'(z)/\Gamma(z)$ denotes the digamma function, and $H_n=\sum_{k=1}^n 1/k$ the $n$th
harmonic number.

{\bf Theorem 7}.  (a) The Ramanujan polynomials satisfy the identities
$$[1-(-1)^n]R_{2s+1}^{(n)}(1)+\sum_{j=1}^{n-1}{n \choose j}\left[{{(n-1)!} \over {(j-1)!}}
-(-1)^j (2s+2)_{n-j}\right]R_{2s+1}^{(j)}(1)$$
$$=(2s+2)_nR_{2s+1}(1),$$
and (b), with $R_{2s}(z)=Q_{2s}(z)$, these generalized Ramanujan polynomials satisfy the
identities
$$[1-(-1)^n]\left[R_{2s}^{(n)}(1)-{1 \over 2^n}R_{2s}^{(n)}\left({1 \over 2}\right)\right]$$
$$+\sum_{j=1}^{n-1}{n \choose j}\left[{{(n-1)!} \over {(j-1)!}}-(-1)^j (2s+2)_{n-j}\right]\left[R_{2s}^{(j)}(1)-{1 \over 2^j}R_{2s}^{(j)}\left({1 \over 2}\right)\right]$$
$$=(2s+2)_n\left[R_{2s}(1)-R_{2s}\left({1 \over 2}\right)\right].$$
(c)  Suppose the following functional equation, as appears in Grosswald's generalization
\cite{grosswald} of Ramanujan's formula pertaining to the Riemann zeta function at odd integer argument
(\cite{gun}, p.\ 945).  For analytic functions $F$ and $S$,
$$F\left(-{1 \over z}\right)-(-1)^\delta \left({z \over i}\right)^r F(z)=S\left({z \over i}\right),$$
where $z$ is in the upper half plane and $r$ is real.
Then
$$[(-1)^n-(-1)^\delta]F^{(n)}(i)-(-1)^\delta (-r)_n i^n F(i)$$
$$+\sum_{j=1}^{n-1}{n \choose j}\left[(-1)^n {{(n-1)!} \over {(j-1)!}}-(-1)^\delta (-1)^{n-j}(-r)_{n-j}
\right]i^{j-n}F^{(j)}(i)={1 \over i^n}S^{(n)}(1).$$
(d) Let $F$ be Zagier's function \cite{zagier} (p.\ 164)
$$F(x)=\sum_{n=1}^\infty{1 \over n}[\psi(nx)-\ln (nx)]=\int_0^\infty \left({1 \over {1-e^{-t}}}-{1 \over t}\right)\ln(1-e^{-xt})dt.$$
Then
$$F'(x)-{1 \over x^2}F\left({1 \over x}\right)=-{\pi^2 \over 6}+{\pi^2 \over {6x^2}}+{{\ln x} \over x},$$
and for $n \geq 2$
$$[1+(-1)^n]F^{(n)}(1)+(-1)^n \sum_{j=1}^{n-1}{n \choose j}{{(n-1)!} \over {(j-1)!}}F^{(j)}(1)
=(-1)^{n-1}(n-1)!\left[{\pi^2 \over 6}n-H_{n-1}\right].$$

The form of the expressions in (a) and (b) shows that no new information is included from the even order
derivative relations, while the form of (d) shows that no new information is included for odd $n$.
Explicitly in terms of Bernoulli summations, the first identity of (a) is the equality of 
$$(s+1)R_{2s+1}(1)=(s+1)\sum_{k=0}^{s+1}{{B_{2k}B_{2s+2-2k}} \over {(2k)!(2s+2-2k)!}}$$
and
$$R_{2s+1}'(1)=\sum_{k=1}^{s+1}{{B_{2k}B_{2s+2-2k}} \over {(2k-1)!(2s+2-2k)!}}.$$

{\bf Theorem 8}.  (General Mellin transform result.)  
Suppose that an analytic function $f$ has an integral representation
$$f(s)={1 \over {\Gamma(s)}}\int_0^\infty t^{s-1}g(t)dt,  ~~~~~~\mbox{Re} ~s>1,$$
for some function $g$.  Putting
$$c_j^M(b)={\sqrt{b} \over {2(2j)!}}\int_0^\infty g(t)[(\sqrt{b}t-1)^{2j}-(\sqrt{b}t+1)^{2j}]dt,$$
we have
$$b^j f(2j)=-c_j^M(b)-\sum_{k=1}^{j-1} {b^k \over {(2j-2k+1)!}} f(2k).$$
Section 3 provides important cases of this result.

In sum, we have generalized the zeta identity of \cite{mcl1}, and the Bernoulli and zeta relations of \cite{mcl2}.  We have noted the relevance of a subset of these results for obtaining summation formulas for double zeta values.  In addition, we give an alternative proof of (1.2) and several other 
generalizations by using integral representations.  Lettington's identity may then be viewed as a
very special case of a recurrence among values of polygamma functions, and beyond this context, 
among values of Dirichlet $L$ series, values of the Lerch zeta function or values of lattice
Dirichlet series.



\medskip
\centerline{\bf Proof of Theorems}
\medskip

{\it Theorem 1}.  (a) We recall the expression for Bernoulli polynomials $B_s(x)=\sum_{k=0}^s {s
\choose k}B_{s-k}x^k$, so that
$$B_j(x)=\sum_{k=0}^j {j \choose k} B_kx^{j-k}.$$
For $n=2j$ even we have the relation $B_n(1/3)=B_n(2/3)=(3^{1-n}-1)B_n/2$ and obtain
$${1 \over 2}(3^{1-2j}-1)B_{2j}=\sum_{k=0}^{2j} {{2j} \choose k}B_k{1 \over 3^{2j-k}}$$
$$={B_0 \over 3^{2j}}+2j{B_1 \over 3^{2j-1}}+\sum_{k=2}^{2j} {{2j} \choose k}B_k{1 \over 3^{2j-k}}$$
$$={1 \over 3^{2j}}-{j \over 3^{2j-1}}+\sum_{m=1}^j {{2j} \choose {2m}}B_{2m}{1 \over 3^{2(j-m)}}$$
$$=-{j \over 3^{2j-1}}+\sum_{m=0}^j {{2j} \choose {2m}}B_{2m}{1 \over 3^{2(j-m)}}.$$ 
The relation (1.1) is then applied. 

(b) is based upon the relation, for $n=2j$ even, $B_n(1/4)=B_n(3/4)=(4^{1-n}-2^{1-n})B_n/2$, so that
$${1 \over 2}(4^{1-2j}-2^{1-2j})B_{2j}=-{j \over 4^{2j-1}}+\sum_{m=0}^j {{2j} \choose {2m}}B_{2m}
{1 \over 4^{2(j-m)}}.$$ 
Then relation (1.1) is again used.

(c) is based upon the relation for $n=2j$ even 
$$B_n\left({1 \over 6}\right)=B_n\left({5 \over 6}\right)={1 \over 2}(6^{1-n}-3^{1-n}-2^{1-n}+1)B_n.$$

(d) is based upon the relation 
$$B_n\left({1 \over {12}}\right)+B_n\left({7 \over {12}}\right)=2^{1-n}B_n\left({1 \over 6}\right)
=2^{-n}(6^{1-n}-3^{1-n}-2^{1-n}+1)B_n.$$

(e) is based upon the relation
$$B_n\left({1 \over {8}}\right)+B_n\left({5 \over {8}}\right)=2^{1-n}B_n\left({1 \over 4}\right)
=2^{-n}(4^{1-n}-2^{1-n})B_n.$$ \qed

{\it Remark}.  The relations that we have employed follow from the symmetry $B_n(1-x)=(-1)^nB_n(x)$
and the multiplication formula $B_n(mx)=m^{n-1}\sum_{k=0}^{m-1} B_n(x+k/m)$.  Hence many more
identities may be developed.

{\it Theorem 2}. (a) We recall the recurrence for Bernoulli numbers
$$B_s=-{1 \over {s+1}}\sum_{k=0}^{s-1} {{s+1} \choose k}B_k.$$
Then
$$\phi_j(2s)={1 \over j^{2s}}\zeta(2s)={{(-1)^{s+1}2^{2s-1}\pi^{2s}} \over {j^{2s}(2s)!}}B_{2s}$$
$$=-{1 \over j^{2s}}{{(-1)^{s+1}2^{2s-1}\pi^{2s}} \over {(2s)!(2s+1)}}\sum_{k=0}^{2s-1}{{2s+1} \choose
k}B_k$$
$$={1 \over j^{2s}}{{(-1)^s(2\pi)^{2s}} \over {2(2s+1)!}}\left[{{2s+1} \choose 1}B_1+\sum_{k=0}^{2s-2}
{{2s+1} \choose k}B_k\right]$$
$$={1 \over j^{2s}}{{(-1)^s(2\pi)^{2s}} \over {2(2s+1)!}}\left[(2s+1)B_1+\sum_{n=0}^{s-1}{{2s+1} \choose {2n}}B_{2n}\right]$$
$$={1 \over j^{2s}}{{(-1)^s(2\pi)^{2s}} \over 2}\left[-{1 \over {2(2s)!}}+\sum_{n=0}^{s-1} {1 \over
{(2s-2n+1)!}}{B_{2n} \over {(2n)!}}\right].$$
The relation (1.1) is used so that
$$\phi_j(2s)={{(-1)^{s+1}2^{2s}} \over {2 j^{2s}}}\left[{\pi^{2s} \over {2(2s)!}}+\sum_{n=1}^s
{{(-1)^{s-n}} \over {(2n+1)!}}{\pi^{2n} \over {2^{2(s-n)-1}}}\zeta(2s-2n)\right]$$
$$={{(-1)^{s+1}2^{2s}} \over {2 j^{2s}}}\left[{\pi^{2s} \over {2(2s)!}}+\sum_{n=1}^s
{{(-1)^{s-n}} \over {(2n+1)!}}{{\pi^{2n}j^{2(s-n)}} \over {2^{2(s-n)-1}}}\phi_j(2s-2n)\right]$$
$$={{(-1)^{s+1}2^{2s}} \over {2 j^{2s}}}\left[{{\pi^{2s}(2s-1)} \over {2(2s+1)!}}+\sum_{n=1}^{s-1}
{{(-1)^{s-n}} \over {(2n+1)!}}{{\pi^{2n}j^{2(s-n)}} \over {2^{2(s-n)-1}}}\phi_j(2s-2n)\right].$$

(b) Theorem 1(a) and relation (1.1) are first used so that
$$-\theta_3(2j){{(2j)!(-1)^{j+1}} \over {(2\pi)^{2j}}}
=-{j \over 3^{2j-1}}+{1 \over 3^{2j}}+\sum_{m=1}^j {{2j} \choose {2m}}{{(-1)^{m+1}(2m)!} \over
{\pi^{2m} 3^{2(j-m)}2^{2m-1}}}\zeta(2m)$$
$$=-{j \over 3^{2j-1}}+{1 \over 3^{2j}}+\sum_{m=0}^{j-1}{{2j} \choose {2m}}{{(-1)^{j-m+1}[2(j-m)]!} \over {\pi^{2(j-m)} 3^{2m}2^{2(j-m)-1}}}\zeta(2j-2m).$$
The relation $\zeta(2s)=(1-3/3^{2s})^{-1}\theta_3(2s)$ is then employed.  

(c) and (d) follow similarly from Theorem 1(b) and 1(c) and relation (1.1). \qed

{\it Theorem 3}.  This follows from Theorem 2(a) and (b) and the result that $\Psi_s$ satisfies the
recurrence (\cite{mcl2}, Lemma 3.1)
$$\Psi_s(\vec{h},\vec{H})=-H_s-\sum_{k=1}^{s-1}h_{s-k}\Psi_k(\vec{h},\vec{H}).$$
\qed

{\it Corollary 1} follows from the multinomial expression for $\zeta(2s)$ of \cite{mcl2}, Lemma 3.3.
\qed

A lower Hessenberg determinant has entries $a_{ij}=0$ for $j-i>1$,
$$A_n=\left|\begin{array}{cccccc}
a_{11} & a_{12}   & 0   & 0 & \ldots & 0 \\
a_{21} & a_{22} & a_{23}   & 0 & \ldots & 0 \\
a_{31} & a_{32} & a_{33} & a_{34} & \ldots & 0 \\
\vdots & \vdots & \vdots & \vdots &\ddots & \vdots \\
a_{n-1,1} & a_{n-1,2} & a_{n-1,3} & a_{n-1,4} & \ldots & a_{n-1,n}\\
a_{n,1} & a_{n,2} & a_{n,3} & a_{n,4} & \ldots & a_{n,n} \\ \end{array}\right|.$$
With $|A_0|=1$ and $|A_1|=a_{11}$, these determinants satisfy the recurrence (e.g., \cite{cahill})
$$\mbox{det}(A_n)=a_{n,n}\mbox{det}(A_{n-1})+\sum_{r=1}^{n-1}\left[(-1)^{n-r}a_{n,r}\mbox{det}
(A_{r-1})\prod_{j=r}^{n-1}a_{j,j+1}\right].$$
The minor corner layered determinants of \cite{mcl1,mcl2} with superdiagonal of all $1$'s and 
comprised of at most $3$ vectors are obviously special cases.  Any time that a recurrence
can be made to take the form just above, a determinant expression may then be written.

The $s \times s$ minor corner layered determinant of \cite{mcl1,mcl2} is given by
$$\Delta_s(\vec{h})=(-1)^s\left|\begin{array}{cccccc}
h_1 & 1   & 0   & 0 & \ldots & 0 \\
h_2 & h_1 & 1   & 0 & \ldots & 0 \\
h_3 & h_2 & h_1 & 1 & \ldots & 0 \\
\vdots & \vdots & \vdots & \vdots &\ddots & \vdots \\
h_{s-1} & h_{s-2} & h_{s-3} & h_{s-4} & \ldots & 1\\
h_s & h_{s-1} & h_{s-2} & h_{s-3} & \ldots & h_1 \\ \end{array}\right|,$$
and we may note that it inverts symmetrically as
$$h_s(\vec{\Delta})=(-1)^s\left|\begin{array}{cccccc}
\Delta_1 & 1   & 0   & 0 & \ldots & 0 \\
\Delta_2 & \Delta_1 & 1   & 0 & \ldots & 0 \\
\Delta_3 & \Delta_2 & \Delta_1 & 1 & \ldots & 0 \\
\vdots & \vdots & \vdots & \vdots &\ddots & \vdots \\
\Delta_{s-1} & \Delta_{s-2} & \Delta_{s-3} & \Delta_{s-4} & \ldots & 1\\
\Delta_s & \Delta_{s-1} & \Delta_{s-2} & \Delta_{s-3} & \ldots & \Delta_1 \\ \end{array}\right|.$$
The recurrence for $\Delta_s(\vec{h})$,
$$\Delta_s(\vec{h})=-\sum_{k=0}^{s-1} h_{s-k}\Delta_k(\vec{h}),$$
is of a standard convolution form.  Indeed, if we let $N(t)=\sum_{j=0}^\infty \Delta_j t^j$,
with $\Delta_0=1$, 
$$d(t)=1+h_1t+h_2t^2+\ldots+h_nt^n\equiv \sum_{j=0}^n h_jt^j,$$
with $h_0=1$, such that $N(t)d(t)=1$, then we have
$$N(t)d(t)=\sum_{m=0}^\infty \sum_{j=0}^m \Delta_jh_{m-j}t^m.$$
This implies the set of equations for $m \geq 0$
$$\sum_{j=0}^m \Delta_jh_{m-j}=\delta_{m0},$$
wherein $\delta_{jk}$ is the Kronecker delta symbol.  Rewriting this equation, we obtain
$$\Delta_m=-\sum_{j=0}^{m-1}\Delta_jh_{m-j}+\delta_{m0}.$$
Thus for $m\geq 1$,
$$\Delta_m=-\sum_{j=0}^{m-1}\Delta_jh_{m-j},$$
just as in \cite{mcl2}, and we next make the connection with the multinomial expansion.

Let $\vec{j}=(j_1,j_2,\ldots,j_n)$, $|\vec{j}|=j_1+j_2+\ldots +j_n$, and
$p(\vec{j})=j_1+2j_2+\ldots +nj_n$.
Then, by the use of geometric series and multinomial expansion,
$${1 \over {d(t)}}=\sum_{k=0}^\infty (-1)^k(h_1t+h_2t^2+\dots+h_nt^n)^k$$
$$=\sum_{k=0}^\infty (-1)^k \sum_{|\vec{j}|=k}{k \choose {j_1,j_2,\ldots,j_k}}
h_1^{j_1}h_2^{j_2}\cdots h_n^{j_n}t^{p(\vec{j})}$$
$$=\sum_{\ell=0}^\infty \sum_{p(\vec{j})=\ell} {{|\vec{j}|} \choose {j_1,j_2,\ldots,j_k}}
(-1)^{|\vec{j}|}h_1^{j_1}h_2^{j_2}\cdots h_n^{j_n}t^\ell.$$
Additionally, we note that the recurrence for $\Delta_m$ may be written in terms of the
companion matrix
$$C_d=\left[\begin{array}{ccccccc}
0 & 0   & 0   & 0 & \ldots & 0 & -h_n \\
1 & 0 & 0   & 0 & \ldots & 0 & -h_{n-1} \\
0 & 1 & 0 & 0 & \ldots & 0 & -h_{n-2} \\
\vdots & \vdots & \vdots & \vdots &\ddots & \vdots & \vdots \\
0 & 0 & 0 & 0 & \ldots & 0 & -h_2\\
0 & 0 & 0 & 0 & \ldots & 1 & -h_1 \\ \end{array}\right].$$

{\it Theorem 4}.  This follows by writing (1.32) of \cite{mcl2} in the form
$${1 \over 2}\left(1-{1 \over 2^{2s}}\right)\zeta(2s)=\sum_{k=1}^{s-1}{1 \over 2^{2(s-k)}}\left(1-
{1 \over 2^{2k}}\right)\zeta(2s-2k)\zeta(2k).$$ \qed

{\it Theorem 5}. We apply the relation
$$E_{2n-1}\left({1 \over 3}\right)=-E_{2n-1}\left({2 \over 3}\right)=-{1 \over {2n}}(1-3^{1-2n}) (2^{2n}-1)B_{2n},$$
so that
$$E_{2n-1}\left({1 \over 3}\right)=-{1 \over {2n}}(1-3^{1-2n}) (2^{2n}-1)B_{2n}$$
$$=\sum_{k=0}^{2n-1} {{2n-1} \choose k}{E_k \over 2^k}\left(-{1 \over 6}\right)^{2n-k-1}$$
$$=\sum_{m=0}^{n-1}{{2n-1} \choose {2m}}{E_{2m} \over {(2m)! }}\left(-{1 \over 6}\right)^{2(n-m)-1}.$$
\qed

{\it Remark}.  Similarly, for instance, we could apply the relation
$$E_{2n}\left({1 \over 6}\right)=E_{2n}\left({5 \over 6}\right)={{(1+3^{-2n})} \over {2^{2n}+1}}E_{2n}.$$

{\it Theorem 6}.  The result follows from relation (1.1) and an identity attributed to Gosper,
$$\sum_{i=0}^n {{(1-2^{1-i})(1-2^{i-n+1})} \over {(n-i)!i!}}B_{n-i}B_i={{(1-n)} \over {n!}}B_n.$$
In order to keep this paper self contained, we provide a proof of this identity.
We use the exponential generating function of Bernoulli polynomials,
$${{te^{xt}} \over {e^t-1}}=\sum_{n \geq 0} {{B_n(x)} \over {n!}}t^n, ~~~~~~|t|<2\pi, \eqno(2.1)$$
so that
$${{te^{t/2}} \over {e^t-1}}=\sum_{n \geq 0} B_n\left({1 \over 2}\right) {t^n \over {n!}},$$
and recall that $B_n(1/2)=(2^{1-n}-1)B_n$.  Then
$${{t^2e^t} \over {(e^t-1)^2}}=\sum_{n=0}^\infty \sum_{m=0}^\infty B_n\left({1 \over 2}\right) 
B_m\left({1 \over 2}\right) {t^{n+m} \over {n!m!}}$$
$$=\sum_{j=0}^\infty \sum_{n=0}^j B_n\left({1 \over 2}\right) B_{j-n}\left({1 \over 2}\right)
{t^j \over {n!(j-n)!}}.$$
We also have
$${{t^2e^t} \over {(e^t-1)^2}}={{2t} \over {e^t-1}}-{d \over {dt}}{t^2 \over {(e^t-1)}}$$
$$=2\sum_{n=0}^\infty B_n {t^n \over {n!}}-{d \over {dt}}\sum_{n=0}^\infty B_n {t^{n+1} \over
{n!}}$$
$$=\sum_{n=0}^\infty [2B_n-(n+1)B_n]{t^n \over {n!}},$$
and we conclude that
$$\sum_{n=0}^j {{B_n\left({1 \over 2}\right) B_{j-n}\left({1 \over 2}\right)} \over
{n!(j-n)!}}=(1-j){B_j \over {j!}}.$$
Gosper's identity then follows. \qed

{\it Remark}.  Using the value $\zeta(0)=-1/2$, the $m=0$ term of the left side of Theorem 6 may
be separated and moved to the right side.  The result is then seen to correspond with a case
presented at the top of p.\ 406 of \cite{basu}.  We have given a distinct method of proof.



{\it Theorem 7}.  The Ramanujan and generalized Ramanujan polynomials respectively satisfy
the functional equations
$$R_{2s+1}(z)=z^{2s+2}R_{2s+1}\left({1 \over z}\right)$$
and
$$R_{2s}(z)-R_{2s}\left({z \over 2}\right)=z^{2s+2}\left[R_{2s}\left({1 \over z}\right)
-R_{2s}\left({1 \over {2z}}\right)\right].$$
We then follow the procedure of \cite{cofpla2002}, noting the relation
$(d/dz)^jz^r=(-1)^j (-r)_jz^{r-j}$, using the product rule for differentiation
and the derivative of a composition of functions.  We arrive at, for instance for part (a),
$${{(-1)^n} \over z^n}\sum_{j=1}^n{n \choose j}{{(n-1)!} \over {(j-1)!}}{1 \over z^j}R_{2s+1}
^{(j)}\left({1 \over z}\right)$$
$$=\sum_{j=0}^n {n \choose j}(-1)^j (2s+2)_jz^{-2s-j-2}R_{2s+1}^{(n-j)}(z).$$
We then set $z=1$, shift the summation index on the right side as $j \to n-j$, and separate
the $j=0$ and $j=n$ terms.  (b) proceeds similarly. 

For (c) we apply the specified functional equation, follow similar steps, and then evaluate at
$z=i$.  

(d) One of the functional equations satisfied by $F$ is \cite{zagier} (p.\ 171)
$$F(x)+F\left({1 \over x}\right)=-{\pi^2 \over 6}x-{\pi^2 \over {6x}}+{1 \over 2}\ln^2 x+C_1,$$
where $C_1=\pi^2/3+2F(1)$ is a constant involving the first Stieltjes constant $\gamma_1$
(e.g., \cite{coffeyst1}), as
$F(1)=-\gamma^2/2-\pi^2/12-\gamma_1$, $\gamma=-\psi(1)$ being the Euler constant.
We again use the derivatives of a composite function and the product rule, and then evaluate at $x=1$.
In the process the relation
\footnote{We may recognize the constants $s(n,2)=(-1)^n(n-1)!H_{n-1}$ as a special case of the
Stirling numbers $s(n,k)$ of the first kind.}
$$\left({d \over {dx}}\right)^j \ln^2 x=2(\ln x-H_{j-1})(-1)^{j-1} {{(j-1)!} \over x^j}$$
is employed.  \qed
\footnote{Other such identities could be developed for Zagier's functions $P(x,y)$, $a(x)$, $A(x,s)$,
and $R(x,y)$.}

{\it Remarks}.  We note that part (c) of the Theorem suffices to write the set of identities implied by (2.3) of \cite{mcl2}, wherein \cite{grosswald}
$$F_s(z)=\sum_{n=1}^\infty \sigma_{-s}(n)e^{2\pi i n z}=\sum_{n=1}^\infty {{\sigma_s(n)} \over n^s}
e^{2\pi inz}=-\zeta(s)-F_s(-z),$$
with $\sigma$ the sum of divisors function.  

$\delta$ takes only the values $0$ and $1$ in \cite{gun}, but that is not necessary in part (c) of the Theorem.

We could write highly related identities at other special points, including $z=-1$ in parts (a) and (b)
and $z=-i$ in part (c).

{\it Theorem 8}.  Is immediate upon applying the integral representation.  \qed

Let $K_s(y)$ be the modified Bessel function of the second kind.
\newline{\bf Corollary 2}.  Define the function for $b\neq 0$ and Re $y>0$
$$c_j^M(b,y)={\sqrt{b} \over {4(2j)!}}\int_0^\infty \exp\left[-{y \over 2}\left(t+{1 \over t}
\right)\right][(\sqrt{b}t-1)^{2j}-(\sqrt{b}t+1)^{2j}]dt.$$
Then
$${b^j \over {(2j-1)!}}K_{2j}(y)=-c_j^M(b,y)-\sum_{k=1}^{j-1} {b^k \over {(2j-2k+1)!}}{{K_{2k}(y)} 
\over {(2k-1)!}}.$$

{\it Proof}.  This follows from the representation for Re $y>0$ and $s \in \mathbb{C}$
$${1 \over {\Gamma(s)}}K_s(y)={1 \over {2\Gamma(s)}}\int_0^\infty \exp\left[-{y \over 2}\left(t+{1 \over t} \right)\right]t^{s-1}dt.$$
\qed

{\bf Corollary 3}.  Define the function for $b\neq 0$ and Re $a>2j$,
$$c_j^M(b,a)={\sqrt{b} \over {2(2j)!}}\int_0^\infty (1+t)^{-a} [(\sqrt{b}t-1)^{2j}-(\sqrt{b}t+1)^{2j}]dt.$$
Then for Re $a>2j$, 
$$b^j\Gamma(a-2j)=-\Gamma(a)c_j^M(b,a)-\sum_{k=1}^{j-1} {{b^k\Gamma(a-2k)} \over {(2j-2k+1)!}}.$$
In fact, this is an identity for a $_3F_2(b)$ hypergeometric function, and we have
$$c_j^M(b,a)=-{b \over {(2j-1)!}}{{\Gamma(a-2)} \over {\Gamma(a)}} {}_3F_2\left(1,{1 \over 2}-j,
1-j;{{3-a} \over 2},2-{a \over 2};b\right).$$

{\it Proof}.  We use a representation for a case of the Beta function $B(x,y)$,
$$\int_0^\infty {t^{x-1} \over {(1+t)^a}}dt={{\Gamma(x)\Gamma(a-x)} \over {\Gamma(a)}}=B(x,a-x),$$
valid for $0<\mbox{Re} ~x<a$, such conditions resulting from maintaining convergence at the 
lower and upper limits of integration, respectively.  We then apply Theorem 8 and multiply 
through by $\Gamma(a)$.

The hypergeometric form for the sum on $k$ may be obtained by shifting the summation index
$k \to k+1$, and then using the duplication formula for Pochhammer symbols, further details being
omitted. \qed

{\it Remarks}.  We may note that a similar result follows for the derivatives with respect to $s$
of the function $K_s(y)/\Gamma(s)$.

For Re $z>0$ and Re $x>0$,
$$x^{-z}={1 \over {\Gamma(z)}}\int_0^\infty e^{-xt}t^{z-1}dt.$$
There follows the identity
$$b^j x^{-2j}=-c_j^M(b,x)-\sum_{k=1}^{j-1} {{b^kx^{-2k}} \over {(2j-2k+1)!}},$$
where
$$c_j^M(b,x)={\sqrt{b} \over {2(2j)!}}\int_0^\infty e^{-xt}[(\sqrt{b}t-1)^{2j}-(\sqrt{b}t+1)^{2j}]dt.$$
Upon evaluating $c_j^M$ using binomial expansion,
$$c_j^M(b,x)=-{1 \over {(2j)!}}\sum_{n=1}^j {{2j} \choose {2n-1}}(2n-1)! {b^n \over x^{2n}},$$
this identity is explicitly verified.  As such, putting $x \to 1/x$ and summing with various
coefficients, this relation becomes a generating identity of, for example, families of polynomials.

\medskip
\centerline{\bf Another proof of the identity (1.2), and other generalizations}
\medskip

In giving a different proof of Lettington's identity (1.2), we make use of a standard integral
representation of the Riemann zeta function,
$$\zeta(s)={1 \over {\Gamma(s)}}\int_0^\infty {t^{s-1} \over {e^t-1}}dt, ~~\mbox{Re} ~s>1, \eqno(3.1)$$
where $\Gamma$ denotes the Gamma function, and Hermite's formula for the Hurwitz zeta function
$$\zeta(s,a)={a^{-s} \over 2}+{a^{1-s} \over {s-1}}+2\int_0^\infty (a^2+y^2)^{-s/2}\sin [s \tan^{-1}(y/a)]{{dy} \over {e^{2\pi y}-1}},$$
holding for all complex $s\neq 1$.

{\bf Lemma 1}.  For integers $j \geq 0$,
$$\int_0^\infty {{(\pi^2+t^2)^j} \over {e^t-1}}\sin \left(2j\tan^{-1} \left({t \over \pi}\right)
\right) dt={{j \pi^{2j+1}} \over {2j+1}}.$$

{\it Proof}.  The $j=0$ case is obvious, and in the following we take $j \geq 1$.  With $v=2\pi y$
and $s=-2j$ in Hermite's formula we obtain
$$\int_0^\infty \left(4a^2+  {v^2 \over \pi^2}\right)^j \sin \left(2j\tan^{-1}{v \over {2\pi a}}
\right) {{dv} \over {e^v-1}}
=\pi 4^j\left[{a^{2j} \over 2}-{a^{2j+1} \over {2j+1}}-\zeta(-2j,a)\right].$$
For $a=1/2$, $\zeta(-2j,1/2)=(2^{-2j}-1)\zeta(-2j)=0$ owing to the trivial zeros of the Riemann
zeta function.  We then find that
$$\int_0^\infty \left(1+  {v^2 \over \pi^2}\right)^j \sin \left(2j\tan^{-1}{v \over \pi}\right)
{{dv} \over {e^v-1}} = {{\pi j}\over {2j+1}},$$
and the Lemma follows.  \qed
\footnote{The reader is invited to provide an alternative proof of Lemma 1 using induction.}

{\it Proof} of (1.2).  We first note the following sum, based upon manipulation of binomial
expansions,
$$\sum_{k=1}^j {{(-1)^k\pi^{2j-2k}t^{2k-1}} \over {(2j-2k+1)!(2k-1)!}}
={\pi^{2j-1} \over {(2j)!}}\sum_{k=1}^j (-1)^k {{2j} \choose {2k-1}}\left({t \over \pi}\right)^{2k-1}$$
$$=-{{i[(\pi-it)^{2j}-(\pi+it)^{2j}]} \over {2\pi (2j)!}}.$$
Also as a prelude, with $y=a \tan^{-1}x$,
$$\sin y = {1 \over {2i}}\left[\left({{1+ix} \over {1-ix}}\right)^{a/2}-\left({{1+ix} \over {1-ix}}\right)^{-a/2}\right].$$
We now use the integral representation (3.1) so that
$$\sum_{k=1}^{j-1} {{(-1)^k \pi^{2j-2k}} \over {(2j-2k+1)!}}\zeta(2k)
=\sum_{k=1}^{j-1} {{(-1)^k \pi^{2j-2k}} \over {(2j-2k+1)!(2k-1)!}}\int_0^\infty {t^{2k-1} \over
{e^t-1}}dt$$
$$=-\int_0^\infty \left[{{2(-1)^j j\pi t^{2j}+t(\pi^2+t^2)^j \sin(2j\tan^{-1}(t/\pi))} \over
{\pi t (2j)!}}\right]{{dt} \over {e^t-1}}.$$
The first term on the right side evaluates to $(-1)^{j-1}\zeta(2j)$ according to (3.1), while
the second term is evaluated by Lemma 1, and (1.2) is shown.  \qed 

{\it Remarks}.  Following steps similar to the proof of Lemma 1, we also find for integers
$j\geq 0$
$$-\int_0^\infty {{(\pi^2+t^2)^{-j}} \over {e^t-1}}\sin \left(2j\tan^{-1} \left({t \over \pi}\right)
\right) dt=\pi\left[{1 \over 2}-{1 \over 2}{1 \over {(1-2j)}}-(1-2^{-2j})\zeta(2j)\right].$$

The method of this section may also be applied to the sums of Theorem 1.

{\bf Theorem 9}.  (a) For Re $a>0$, define
$$c_j(a)={1 \over {\pi(2j)!}}\int_0^\infty (\pi^2+t^2)^j \sin(2j\tan^{-1}(t/\pi))
{{e^{-(a-1)t}dt} \over {e^t-1}}.$$
Then
$$\zeta(2j,a)=(-1)^{j-1}\left[c_j(a)+\sum_{k=1}^{j-1}{{(-1)^k\pi^{2j-2k}} \over
{(2j-2k+1)!}}\zeta(2k,a)\right].$$
(b) Define for Re $a>0$ and $b \neq 0$
$$c_j(a,b)={\sqrt{b} \over {\pi(2j)!}}\int_0^\infty (\pi^2+bt^2)^j \sin\left(2j\tan^{-1}\left({{\sqrt{b}t} \over \pi}\right)\right)
{{e^{-(a-1)t}dt} \over {e^t-1}}.$$
Then
$$b^j\zeta(2j,a)=(-1)^{j-1}\left[c_j(a,b)+\sum_{k=1}^{j-1}{{(-1)^k\pi^{2j-2k}b^k} \over
{(2j-2k+1)!}}\zeta(2k,a)\right].$$

{\it Proof}.  We use the integral representation for Re $a>0$  
$$\zeta(s,a)={1 \over {\Gamma(s)}}\int_0^\infty {{t^{s-1}e^{-(a-1)t}} \over {e^t-1}}dt, ~~\mbox{Re} ~s>1,$$
and follow steps similar to the above in reproving (1.2). \qed

{\it Remark}.  We suspect that the constants $c_j(a=1/2)$ may be evaluated explicitly.

We now introduce Dirichlet $L$-functions $L(s,\chi)$ (e.g., \cite{ireland}, Ch. 16), that are 
known to be expressible as linear combinations of Hurwitz zeta functions.  
We let $\chi_k$ be a Dirichlet character modulo $k$, 
and have
$$L(s,\chi)=\sum_{n=1}^\infty {{\chi_k(n)} \over n^s}, ~~~~\mbox{Re} ~s >1. $$
This equation holds for at least Re $s>1$.  If $\chi_k$ is a nonprincipal character, 
then convergence obtains for Re $s>0$.  The $L$ functions are extendable to the whole complex
plane, satisfy functional equations, and have integral representations.


{\bf Corollary 4}.  Define the constants
$$b_j^{(k)}={1 \over k^s}\sum_{m=1}^k \chi_k(m)c_j\left({m \over k}\right).$$
Then for $L$ functions of characters modulo $k$,
$$L(2j,\chi)=(-1)^{j-1}\left[b_j^{(k)}+\sum_{\ell=1}^{j-1}{{(-1)^k\pi^{2j-2\ell}} \over
{(2j-2\ell+1)!}}L(2\ell,\chi)\right].$$

{\it Proof}.  This follows from Theorem 7(a) and the relation
$$L(s,\chi)={1 \over k^s}\sum_{m=1}^k \chi_k(m) \zeta\left(s,{m \over k}\right).$$
\qed

We may now emphasize the Lettington relation (1.2) as a recurrence among special values of polygamma
functions $\psi^{(j)}$.  
\footnote{We recall the definitions $\psi(z)=\Gamma'(z)/\Gamma(z)$ and
$\psi^{(j)}(z)=(d^j/dz^j)\psi(z)$, and the functional equation $\psi^{(j)}(x+1)=(-1)^jj!/x^{j+1}+
\psi^{(j)}(x)$.}
For we have the relation
$$\zeta(2j,a)={1 \over {(2j-1)!}}\psi^{2j-1}(a),$$
and the example evaluations for $j>0$
$$\psi^{(j)}(1)=(-1)^j (j-1)!\zeta(j)$$
and
$$\psi^{(j)}\left({1\over 2}\right)=(-1)^j (j-1)!(2^j-1)\zeta(j).$$ 

The Lerch zeta function is given by (e.g., \cite{grad} p.\ 1075)
$$\Phi(z,s,a)=\sum_{n=0}^\infty {z^n \over {(n+a)^s}}.$$
This series converges for $a \in \mathbb{C}$ not a negative integer and all $s \in \mathbb{C}$ when $|z|<1$, and for Re $s>1$ when $|z|=1$
and has the integral representation generalizing Hermite's formula
$$\Phi(z,s,a)={a^{-s} \over 2}+\int_0^\infty {z^t \over {(t+a)^s}}dt
-2\int_0^\infty (a^2+y^2)^{-s/2}\sin [y\ln z-s \tan^{-1}(y/a)]{{dy} \over {e^{2\pi y}-1}},$$
for Re $a>0$.  It satisfies the functional equation
$$\Phi(z,s,a)=z^n\Phi(z,s,n+a)+\sum_{k=0}^{n-1}{z^k \over {(k+a)^s}}.$$

{\bf Theorem 10}. Define for Re $a>0$ and $b \neq 0$
$$c_j(a,b,z)={\sqrt{b} \over {\pi(2j)!}}\int_0^\infty (\pi^2+bt^2)^j \sin\left(2j\tan^{-1}\left({{\sqrt{b}t} \over \pi}\right)\right)
{{e^{-(a-1)t}dt} \over {e^t-z}}.$$
Then
$$b^j\Phi(z,2j,a)=(-1)^{j-1}\left[c_j(a,b,z)+\sum_{k=1}^{j-1}{{(-1)^k\pi^{2j-2k}b^k} \over
{(2j-2k+1)!}}\Phi(z,2k,a)\right].$$

{\it Proof}.  We use the integral representation for Re $a>0$  
$$\Phi(z,s,a)={1 \over {\Gamma(s)}}\int_0^\infty {{t^{s-1}e^{-(a-1)t}} \over {e^t-z}}dt, ~~\mbox{Re} ~s>1.$$
\qed

We let $_pF_q$ be the generalized hypergeometric function and
$$\mbox{Li}_s(z)=\sum_{k=1}^\infty {z^k \over k^s}=z\Phi(z,s,1),$$
for $s \in \mathbb{C}$ and $|z|<1$ or Re $s>1$ and $|z|=1$ be the polylogarithm function.
\newline{\bf Corollary 5}.  (a) For Re $a>0$ and $b \neq 0$,
$$b^ja^{-2j} {}_{2j+1}F_{2j}(1,a,\ldots,a;a+1,\ldots,a+1;z)
=(-1)^{j-1}\left[c_j(a,b,z)\right.$$
$$\left.+\sum_{k=1}^{j-1}{{(-1)^k\pi^{2j-2k}b^k} \over {(2j-2k+1)!}}a^{-2k}
{}_{2k+1}F_{2k}(1,a,\ldots,a;a+1,\ldots,a+1;z)\right],$$
and (b) for $|z| \leq 1$,
$$b^j\mbox{Li}_{2j}(z)=(-1)^{j-1}\left[zc_j(1,b,z)+\sum_{k=1}^{j-1}{{(-1)^k\pi^{2j-2k}b^k} \over {(2j-2k+1)!}}\mbox{Li}_{2k}(z)\right].$$

{\it Proof}.  (a) follows from the relation for integers $k\geq 1$
$$\Phi(z,k,a)=a^{-k} {}_{k+1}F_k(1,a,\ldots,a;a+1,\ldots,a+1;z)$$
and (b) follows from $\Phi(z,2j,1)=\mbox{Li}_{2j}(z)/z$.  \qed

Other functions will have a recurrence similar to, and in several cases again generalizing, (1.2),
and we conclude this section noting that this applies to Eisenstein series.
Let
$$f_1(\tau,t)={{\cosh^2(\tau t/2)} \over {1-2e^{-t}\cosh(\tau t)+e^{-2t}}},$$
$$f_2(\tau,t)={{\cos^2(t/2)} \over {1-2e^{i\tau t}\cos t+e^{2i\tau t}}},$$
and
$$\tilde{E}_s=\lim_{K \to \infty}\sum_{|m|,|n|\leq K}{1 \over {(m\mu+n \nu)^s}},$$
where $\mu,\nu \in \mathbb{C}$ and $\tau=\nu/\mu \notin \mathbb{R}$.
Then there  is a recurrence for values $\tilde{E}_{2j}$, owing to the integral 
representation (\cite{dienstfrey}, Theorem 5, $\mu=1$, $\nu=\tau$)
$$\tilde{E}_s(\tau)=\cos\left({\pi \over 2}s\right){4 \over {\Gamma(s)}}\int_0^\infty t^{s-1}
[e^{-is\pi/2}e^{-t}f_1(\tau,t)+e^{i\tau t}f_2(\tau,t)]dt, ~~~~\mbox{Re} ~s>2.$$
Here the ratio $\tau$ is in the fundamental region $-1/2<\mbox{Re}~\tau\leq 1/2$,
Im $\tau>0$, $|\tau| \geq 1$, and if $|\tau|=1$, then Re $\tau \geq 0$.  By Corollary 6 of
\cite{dienstfrey}, the summatory conditionally convergent case of $s=2$ is also given by this
integral.  Walker \cite{walker} found a remarkable formula for $\tilde{E}_2(\tau)$ in terms of the Dedekind $\eta(\tau)=e^{i\pi \tau/12}\prod_{n=1}^\infty (1-e^{2\pi i \tau n})$ function, when
summing over increasing disks.

Omitting further details, we arrive at:
\newline{\bf Theorem 11}.  Define the function, for $\tau$ in the fundamental region in the
upper half plane, and $b \neq 0$,
$$C_j^E(b,\tau)={{4\sqrt{b}} \over {\pi(2j)!}}\int_0^\infty \left[(\pi^2+bt^2)^j \sin\left(2j
\tan^{-1}\left({{\sqrt{b}t} \over \pi}\right)\right)e^{-t}f_1(\tau,t)\right.$$
$$+\left.i(\pi^2-bt^2)^j \sin\left(2j \tan^{-1}\left({{\sqrt{-b}t} \over \pi}\right)\right)e^{it\tau}f_2(\tau,t)\right]dt.$$
Then for $\tau$ in the fundamental region,
$$b^j\tilde{E}_{2j}(\tau)=(-1)^{j-1}\left[C_j^E(b,\tau)+\sum_{k=1}^{j-1}{{(-1)^k\pi^{2j-2k}b^k} \over
{(2j-2k+1)!}}\tilde{E}_{2k}(\tau)\right].$$

{\it Remarks}.  It appears that a generalized Hermite formula for $\tilde{E}_s$ should exist,
and we may guess that it contains a half-line integral term something like
$$\int_0^\infty [(\pi^2+t^2)^{-s/2}\sin(s \tan^{-1}(t/\pi))e^{-t}f_1(\tau,t)
+ i(\pi^2-t^2)^{-s/2}\sin(s \tan^{-1}(t/\pi))e^{it\tau}f_2(\tau,t)]dt.$$
In fact, such a formula should follow from the contour integral representation of Theorem 7 of
\cite{dienstfrey}, and as a precursor we have the following result for the Hurwitz zeta function.
Hermite-type formulas are usually found via Plana summation.  However, the following shows that
this is not necessary.

{\bf Theorem 12}.  For Re $s>1$, Re $a>0$, and $0<c<1$,
$$\zeta(s,a)=-\int_0^\infty \left[\cos\left(s\tan^{-1}\left({t \over {c+a-1}}\right)\right)
\cos \pi c \sin \pi c \right.$$
$$\left.-\sinh \pi t \cosh \pi t\sin \left(s\tan^{-1}\left({t \over {c+a-1}}\right)\right)\right]$$
$$\times \left[(c+a-1)^2+t^2]^{s/2}\left(\cosh^2 \pi t-\cos^2 \pi t\right)\right]^{-1}dt.$$

{\it Proof}.  For Re $s>1$, we have
$$\zeta(s,a)=\pi \sum_{k=0}^\infty \mbox{Res} \left.(t+a)^{-s} \cot \pi t \right|_{t=k},$$
leading to
$$\zeta(s,a)={i \over 2}\int_{c-i\infty}^{c+i\infty} (t+a-1)^{-s} \cot \pi t ~dt$$
$$=-{1 \over 2}\int_{-\infty}^\infty (c+a-1+it)^{-s} \cot[\pi(c+it)]dt.$$
The contributions for negative and positive $t$ are then combined, using for instance
$$\cot \pi(c-it)={{\cos \pi c\cosh \pi t-i\sin \pi c \sinh \pi t} \over {\cosh \pi t\sin \pi c+i
\cos \pi c \sinh \pi t}}.$$
\qed

We may similarly consider certain lattice Dirichlet series (Kronecker series) 
$$G(s,\chi)=\sum_{\omega \in \Lambda \backslash \{0\}} {{\chi(\omega)} \over {|\omega|^{2s}}}, ~~~~~~
\chi(\omega)=e^{i(m\mu \alpha+n \nu \beta)},$$
with $\alpha$ and $\beta$ real and $\Lambda \in \mathbb{C}$ a lattice.  Letting
$$|\omega_{m,n}|^2=|m\mu+n\nu|^2=Q(m,n),$$
there is the integral representation \cite{dienstfrey}
$$G(s,\chi)={1 \over {\Gamma(s)}}\int_0^\infty t^{s-1} \sum_{\omega \in \Lambda\backslash \{0\}} \chi(\omega)e^{-tQ(m,n)}dt,$$
leading to the following result, whose proof is omitted.

{\bf Theorem 13}.  Define for $b \neq 0$
$$C_j^G(b,\chi)={\sqrt{b} \over {\pi(2j)!}}\sum_{\omega \in \Lambda \backslash \{0\}}\chi(\omega)
\int_0^\infty (\pi^2+bt^2)^j \sin\left(2j\tan^{-1}\left({{\sqrt{b}t} \over \pi}\right)\right)
e^{-tQ(m,n)}dt.$$
Then
$$b^jG(2j,\chi)=(-1)^{j-1}\left[C_j^G(b,\chi)+\sum_{k=1}^{j-1}{{(-1)^k\pi^{2j-2k}b^k} \over
{(2j-2k+1)!}}G(2k,\chi)\right].$$

Let $\sigma_k(n)$ be the sum of divisors function, the sum of the powers $d^k$ of the positive divisors of
$n$, i.e., $\sigma_k(n)=\sum_{d|n} d^k$.  With nome $q=e^{i\pi \tau}$, there are the series
representations
$$\tilde{E}_{2k}(\tau)=2\zeta(2k)+{{2(2\pi i)^{2k}} \over {(2k-1)!}}\sum_{n=1}^\infty \sigma_{2k-1}(n)
q^n$$
$$=2\zeta(2k)+{{2(2\pi i)^{2k}} \over {(2k-1)!}}\sum_{n=1}^\infty {{n^{2k-1} q^n} \over {1-q^n}}.$$
These series provide another means by which to establish the equivalent of Theorem 11, and of the
various generalizations of Theorem 1 for $\tilde{E}_{2k}(\tau)$.

\medskip
\centerline{\bf Integral and series representations for the Hurwitz numbers $\tilde{H}_n$}
\medskip

The Hurwitz numbers $\tilde{H}_n$ \cite{carlitz62,hurwitz} occur in the Laurent expansion of a certain Weierstrass $\wp$ function about the origin, and are highly analogous to the Bernoulli numbers.
The $\wp$ function in question has periods $\tilde{\omega}$ and $\tilde{\omega}i$, where $\tilde{\omega}$
is the Beta function value
$$\tilde{\omega}=2\int_0^1 {{dx} \over \sqrt{1-x^4}}={1 \over 2}B\left({1 \over 4},{1 \over 2}\right)
={\sqrt{\pi} \over 2}{{\Gamma(1/4)} \over {\Gamma(3/4)}},$$
and satisfies the nonlinear differential equations
$$\wp'(z)^2=4\wp(z)^3-4\wp(z), ~~~~~~~~\wp''(z)=6\wp(z)^2-2.$$
Then this $\wp$ function expands about $z=0$ as
$$\wp(z)={1 \over z^2}+\sum_{n=2}^\infty {{2^n \tilde{H}_n} \over n}{z^{n-2} \over {(n-2)!}}. \eqno(4.1)$$
From the property $\wp(iz)=-\wp(z)$, and the evenness of $\wp(z)$, it follows that $\tilde{H}_n=0$ unless $n$ is a multiple of $4$,
the first few values being $\tilde{H}_4=1/10$, $\tilde{H}_8=3/10$, and $\tilde{H}_{12}=567/130$.
The analogous expansion for Bernoulli numbers is
$${1 \over {\sin^2 x}}={1 \over x^2}+\sum_{n=2}^\infty {{(-1)^{n/2-1}2^n B_n} \over n}{x^{n-2} \over {(n-2)!}}. \eqno(4.2)$$

The analogy between $B_{2n}$ and $\tilde{H}_{4n}$ is furthered by comparing the Riemann zeta function
values
$$\sum_{r \in \mathbb{Z}^+}{1 \over r^{2n}}=(-1)^{n-1}{{(2\pi)^{2n}} \over {(2n)!}}B_{2n}, ~~n \geq 1,$$
with the sum over Gaussian integers
$$\sum_{\lambda \in \mathbb{Z}+i\mathbb{Z}\backslash(0,0)}{1 \over \lambda^{4n}}={{(2\tilde{\omega})^{4n}} \over {(4n)!}}\tilde{H}_{4n}.$$         
In addition, Hurwitz \cite{hurwitz} proved a Clausen-von Staudt type theorem for $\tilde{H}_n$.

An integral representation for $\tilde{H}_n$ would have many applications, including the
development of various recurrences, again in analogy to the previous sections of this paper.
Therefore, we present the following.
{\newline \bf Theorem 14}.  For $k\geq 1$ and the function $f_1$ defined as in the previous section,
$$\tilde{H}_{2k+2}={{2(k+1)} \over {4^k \tilde{\omega}^{2k+2}}}[1+(-1)^{k+1}]\int_0^\infty e^{-t}f_1(i,t)t^{2k+1}dt.$$

{\it Proof}.  We make use of Theorem 10 of \cite{dienstfrey}, such that
\footnote{Note that in (33) and (34) of \cite{dienstfrey}, the $\tau$ and $\lambda$ arguments of
the functions $f_1$ and $f_2$ need to be reversed.}
$$\wp(z,\tau)={1 \over z^2}+8\int_0^\infty t\left[e^{-t}\sinh^2\left({{zt} \over 2}\right)f_1(\tau,t)
+e^{i t \tau}\sin^2\left({{zt} \over 2}\right)f_2(\tau,t)\right]dt,$$
with $\tau$ in the fundamental region and $z$ in a domain containing the origin.  With the
Maclaurin series
$$\sin^2 x={1 \over 2}(1-\cos 2x)=\sum_{k=1}^\infty (-1)^{k+1} {{2^{2k-1}x^{2k}} \over {(2k)!}},$$
and
$$\sinh^2 x={1 \over 2}(\cosh 2x-1)=\sum_{k=1}^\infty {{2^{2k-1}x^{2k}} \over {(2k)!}},$$
we may write
$$\wp(z)-{1 \over z^2}=4\sum_{k=1}^\infty{z^{2k} \over {(2k)!}}\int_0^\infty [e^{-t}f_1(\tau,t)+e^{it \tau}(-1)^{k+1}f_2(\tau,t)]t^{2k+1}dt. \eqno(4.3)$$
For the square lattice with $\tau=i$, $f_1(i,t)=f_2(i,t)$. 
We then equate like powers of $z$ with those of the expansion (4.1) and apply the scaling
$\wp(x|\mu,\nu)=\wp(x/\mu|1,\tau)/\mu^2$ \cite{dienstfrey} (p.\ 145). \qed 

{\bf Corollary 6}.  Define the constants
$$c_{\ell,j}^{(k)}=\int_0^\infty e^{-\ell t}t^{4k-1}\cos^{\ell-(2j+1)}t ~dt.$$
Then
$$\tilde{H}_{4k}={k \over 4^{2(k-1)}}{1 \over \tilde{\omega}^{4k}}\sum_{\ell=1}^\infty \sum_{j=0}^
{(\ell-1)/2} 2^{\ell-(2j+1)}{{\ell-j-1} \choose j}(-1)^j(c_{\ell,j}^{(k)}+c_{\ell+1,j}^{(k)}).$$

{\it Proof}.  In the integral representation for $\tilde{H}_{4n}$, we may note the simplified product
$$e^{-t}f_1(i,t)={{\cos t+1} \over {4(\cosh t-\cos t)}}.$$
From the identity
$$\sum_{k=1}^\infty e^{-kt}\sin kx={{\sin x} \over {2(\cosh t-\cos x)}},$$
we have
$$\tilde{H}_{4k}={k \over 4^{2(k-1)}}{1 \over \tilde{\omega}^{4k}}\sum_{\ell=1}^\infty \int_0^\infty
{{(\cos t+1)} \over {\sin t}}e^{-\ell t}t^{4k-1}\sin \ell t ~dt.$$
We then apply the identity
$$\sin nx =\sin x\sum_{j=0}^{(n-1)/2} (-1)^j {n \choose {2j+1}}\cos^{n-(2j+1)}x \sin^{2j}x$$
$$=\sin x\sum_{j=0}^{(n-1)/2} (-1)^j 2^{n-(2j+1)}{{n-j-1} \choose j}\cos^{n-(2j+1)}x.$$
\qed   

{\it Remarks}.  The form of Theorem 14 verifies that $\tilde{H}_{4j+2}=0$ and $\tilde{H}_{4n}\neq 0$.
The recurrence relation for $\tilde{H}_{4n}$ obtained from the differential equation for $\wp(z)$ for the square lattice is 
$$(2n-3)(4n-1)(4n+1)\tilde{H}_{4n}=3\sum_{j=1}^{n-1}(4j-1)(4n-4j-1){{4n} \choose {4j}}\tilde{H}_{4j}
\tilde{H}_{4(n-j)}.$$

By rewriting the result of Theorem 14, we have
$$\tilde{H}_{4k}={{8k} \over {4^{2k-1}\tilde{\omega}^{4k}}}\int_0^\infty e^{-t}f_1(i,t)t^{4k-1}dt,$$
and this highly suggests the study of sums of the form
$$\sum_{k=1}^{n-1} {{4n+p} \choose {4k}}{4^{2k-1} \over {8k}}\tilde{H}_{4k},$$
with $p=0,1,2,3$. 
Indeed, we obtain the following.
\newline{\bf Theorem 15}.  Define the constants
$$C_n^H=\int_0^\infty {1 \over t}e^{-t}f_1(i,t)\left[-1+{}_4F_3\left({1 \over 4}-n,{1 \over 2}
-n,{3 \over 4}-n,-n;{1 \over 4},{1 \over 2},{3 \over 4};{t^4 \over \tilde{\omega}^4}\right)\right]dt.$$
Then
$${4^{2n-1} \over {8n}}\tilde{H}_{4n}=C_n^H-\sum_{k=1}^{n-1} {{4n} \choose {4k}}{4^{2k-1} \over {8k}}\tilde{H}_{4k}.$$

{\it Proof}.  We apply the integral representation of the previous Theorem, 
$$\sum_{k=1}^{n-1} {{4n} \choose {4k}}{4^{2k-1} \over {8k}}\tilde{H}_{4k}=\sum_{k=1}^{n-1}
{{4n} \choose {4k}}{1 \over \tilde{\omega}^{4k}}\int_0^\infty e^{-t}f_1(i,t)t^{4k-1}dt,$$
with
$${{4n} \choose {4k}}={{(-4n)_{4k}} \over {(4k)!}}.$$
In order to reach hypergeometric form, we then apply the quadriplication formulas for the Pochhammer symbol and Gamma function, so that
$$(-4n)_{4k}=4^{4k}(-n)_k\left(-n+{1 \over 4}\right)_k\left(-n+{1 \over 2}\right)_k
\left(-n+{3 \over 4}\right)_k,$$
and
$$(4k)!=\Gamma(4k+1)
=(2\pi)^{-3/2}4^{4k+1/2}\Gamma\left(k+{1 \over 4}\right)\Gamma\left(k+{1 \over 2}\right)\Gamma\left(k+{3 \over 4}\right)k!$$
$$=4^{4k}\left({1 \over 4}\right)_k\left({1 \over 2}\right)_k\left({3 \over 4}\right)_k k!.$$
\qed

Another very useful application of an integral representation can be in developing asymptotic 
relations.  Indeed, we have the following.
{\newline \bf Corollary 7}.  As $n \to \infty$,
$$\tilde{H}_{4n} \sim 4 {{(4n)!} \over {(2\tilde{\omega})^{4n}}},$$
with the refinement
$$\tilde{H}_{4n} \sim 4 {{(4n)!} \over {2^{6n}\tilde{\omega}^{4n}}}[(-1)^k+2^{2n}].$$

{\it Proof}.  We have
$$\tilde{H}_{4k} \sim {{8k} \over {4^{2k-1}\tilde{\omega}^{4k}}}\int_0^\infty e^{-t} \cos^2(t/2)t^{4k-1}
dt$$
$$={{4k} \over {4^{2k-1}\tilde{\omega}^{4k}}}\int_0^\infty e^{-t} (\cos t+1)t^{4k-1}dt$$
$$={{8k} \over {4^{2k-1}\tilde{\omega}^{4k}}}{1 \over 4^{2k+1}}[(1-i)^{4k}+(1+i)^{4k}+2^{4k+1}]
\Gamma(4k)$$
$$={{(4k)!} \over {4^{2k-1}2^{2k}\tilde{\omega}^{4k}}}[(-1)^k+2^{2k}].$$
\qed

The following collects various integral representations for lattice sums $S_r(\Lambda)=\sum_{\omega \in
\Lambda\backslash(0,0)} 1/\omega^r$ and the $g_2$ and $g_3$ invariants of the Weierstrass $\wp$
function.  
{\newline \bf Corollary 8}. Let the functions $f_1$ and $f_2$ be defined as in the previous section.
Then (a)
$$S_{2k}(\Lambda)={4 \over {(2k+1)!}}\int_0^\infty[e^{-t}f_1(\tau,t)+(-1)^{k+1}e^{i\tau t}f_2(\tau,t)]t^{2k+1}dt,$$
(b) 
$$g_2={{40} \over \mu^4}\int_0^\infty [e^{-t}f_1(\tau,t)+e^{i\tau t}f_2(\tau,t)]t^3dt,$$
and
(c)
$$g_3={{14} \over {3\mu^6}}\int_0^\infty [e^{-t}f_1(\tau,t)-e^{i\tau t}f_2(\tau,t)]t^5dt.$$
The discriminant of the polynomial $4x^3-g_2x-g_3$ is given by
$\Delta(\Lambda)=g_2(\Lambda)^3-27g_3(\Lambda)^2$ and the $j$-invariant of the lattice $\Lambda$
is $j(\Lambda)=1728g_2(\Lambda)^3/\Delta(\Lambda)$.

{\it Proof}.  (a) follows from the work in proving Theorem 14.  We have
$S_4=g_2/60$ and $S_6=g_3/140$ (e.g., \cite{armitage}, p. 166), so that
$\wp(z,\Lambda)-1/z^2=g_2z^2/20+g_3z^4/28+O(z^6)$.  Then for (b),
$$\left(\wp(z,\Lambda)-{1 \over z^2}\right)''_{z=0}={g_2 \over {10}},$$
while for (c),
$$\left(\wp(z,\Lambda)-{1 \over z^2}\right)^{(iv)}_{z=0}={{6g_3} \over 7}.$$
\qed

{\it Remarks}.  For the square lattice with $\tau=i$ and $\mu=\tilde{\omega}$ (and the scaling of \cite{dienstfrey}) we obtain as expected $g_2=4$ and $g_3=0$.

If $F$ is a field of characteristic different from $2$ or $3$, then elliptic curves $E(F)$ over $F$ are
isomorphic if they have the same $j$-invariant.  The $j$-invariant is a modular function of weight $0$,
and it can be shown that any modular function of weight $0$ must be a rational function of $j$.

We also obtain integral representations for the coefficients of the Laurent expansion about $z=0$ for the Weierstrass zeta function $\zeta_w$ from Corollary 8 since $\zeta_w(z,\tau)=-\int^z \wp(s,\tau)ds =1/z-S_4z^3-S_6z^5 +O(z^7)$.

Our results are also relevant to series expansions of the solutions of the first Painlev\'{e} equation
${{d^2u} \over {dz^2}}=6u^2+z$.  This is because Boutroux \cite{boutroux} showed that, for large $|z|$, with the scaled variables $U = z^{-1/2}u$ and $Z ={4\over 5}z^{5/4}$, the solution of this equation behaves asymptotically like the Weierstrass function, $U\sim \wp$, which satisfies the second order differential equation $\wp''=6\wp^2-g_2/2$.

Absent a (much) more general future result, the following illustrates other identities satisfied by
the Hurwitz numbers.
\newline{\bf Proposition 1}.
$$-30\tilde{H}_4^2+\tilde{H}_8=0, ~~~~~~
6\tilde{H}_4^3-{9 \over 5}\tilde{H}_4\tilde{H}_8+{{52} \over {4725}}\tilde{H}_{12}=0,$$
and
$$68\tilde{H}_4^2\tilde{H}_8-{{16} \over 5}\tilde{H}_8^2-{{1408} \over {945}}\tilde{H}_4\tilde{H}_{12}
+{{901} \over {315315}}\tilde{H}_{16}=0.$$

{\it Proof}.  We write the duplication formula of the $\wp$ function in the form
$$\wp(2z)[\wp'(z)]^2={1 \over 4}[\wp''(z)]^2-2[\wp'(z)]^2\wp(z),$$
wherein polar terms of $O(1/z^8)$ and $O(1/z^4)$ cancel.  Then comparing coefficients of the terms $O(1)$, $O(z^4)$, and $O(z^8)$ gives the result.  \qed

{\bf Theorem 16}. (Series representation of the Hurwitz numbers).  For $m \geq 1$,
$$\tilde{H}_{4m}=\left({\pi \over \tilde{\omega}}\right)^{4m}\left[-B_{4m}+8m \sum_{n=1}^\infty
{n^{4m-1} \over {e^{2\pi n}-1}}\right].$$

{\it Proof}.  The $\wp$ function has the Fourier expansion \cite{basoco}
$$\wp(z,\tau)=-2\left({1 \over 6}+\sum_{n=1}^\infty {1 \over {\sin^2(n\pi \tau)}}\right)
+{\pi^2 \over {\sin^2 \pi z}}
-8\pi^2 \sum_{n=1}^\infty {{n \cos(2\pi n z)} \over {e^{-2\pi i \tau n}-1}}.$$
This expression is expanded in powers of $z$ with $\tau=i$, appropriately scaling in terms of
$\tilde{\omega}$, and the use of (4.2) and the following.
\newline{\bf Lemma 2}.  
$$\sum_{n=1}^\infty {1 \over {\sin^2(n\pi i)}}={1 \over \pi^2}\sum_{n=1}^\infty \Gamma^2(in)\Gamma^2(1-in)
=-{1 \over 6}+{1 \over {2\pi}}.$$
{\it Proof of Lemma 2}.  
There is the expansion
$$\csc^2 \pi z={1 \over {\sin^2 \pi z}}={1 \over {\pi^2 z^2}}+{2 \over \pi^2}\sum_{k=1}^\infty
{{z^2+k^2} \over {(z^2-k^2)^2}}.$$
Then
$$\sum_{n=1}^\infty {1 \over {\sin^2(n\pi i)}}=-{1 \over \pi^2}\sum_{n=1}^\infty {1 \over n^2}
+{2 \over \pi^2}\sum_{n=1}^\infty \sum_{k=1}^\infty {{(k^2-n^2)} \over {(k^2+n^2)^2}}$$
$$=-{1 \over 6}+{2 \over \pi^2}\sum_{k=1}^\infty \left[-{1 \over {2k^2}}+{\pi^2 \over 2}
\mbox{csch}^2(\pi k)\right].$$ \qed

There follows 
$$\sum_{m=1}^\infty {{2^{4m} \tilde{H}_{4m}} \over {4m(4m-2)!}}z^{4m-2}={1 \over \tilde{\omega}^2}
\left[-\pi +2B_2 \pi^2 -\sum_{m=1}^\infty {{2^{4m} B_{4m}\pi^{4m}} \over {4m(4m-2)!}}{z^{4m-2}
\over {\tilde{\omega}^{4m-2}}}\right.$$
$$\left. -8\pi^2\sum_{j=0}^\infty {{(-1)^j} \over {(2j)!}} {{(2\pi)^{2j}} \over \tilde{\omega}^{2j}}
\sum_{n=1}^\infty {n^{2j+1} \over {(e^{2\pi n}-1)}}z^{2j}\right].$$
Equating the coefficients of $z^{4m-2}$ gives the Theorem.  \qed

\pagebreak
\centerline{\bf The equianharmonic case for $\wp$}
\medskip

We have covered the lemniscatic case of the $\wp$ function, for which the lattice is a 
certain square.  We are concerned in this section with the equianharmonic case, for which
the lattice is composed of equilateral triangles.  
\footnote{Further background on these two cases may be found in section 7.5 of \cite{walkerbk}.}
So we now consider a $\wp$ function for
which $g_2=0$, $g_3=1$, with half periods
$$\omega_1=\int_{4^{1/3}}^\infty {{dx} \over \sqrt{4x^3-1}}={{\Gamma^3(1/3)} \over {4\pi}},$$
and $\omega_2=e^{i\pi/3}\omega_1$.  Then $\tau=\omega_2/\omega_1=e^{i\pi/3}$ is a $6$th root
of unity.

We write the expansion
$$\wp(z;0,1)={1 \over z^2}+\sum_{n=1}^\infty {{(6n-1)} \over {(2\omega_1)^{6n}}}S_{6n}(e^{i\pi/3})
z^{6n-2}, \eqno(5.1)$$
such that $S_{2n}(e^{i\pi/3})=0$ unless $n \equiv 0$ mod 3.

{\bf Theorem 17}. Let the functions $f_1$ and $f_2$ be as given in section 3.  Then (a)
$$S_{6n}(e^{i\pi/3})={4 \over {(6n-1)!}}\int_0^\infty [e^{-t}f_1(e^{i\pi/3},t)+e^{it\exp(i\pi/3)}
(-1)^n f_2(e^{i\pi/3},t)]t^{6n-1}dt,$$
and (b)
$$S_{6n}(e^{i\pi/3})\sim 6+2\left(-{1 \over {27}}\right)^n.$$

{\it Proof}. Part (a) follows by using the intermediate result (4.3) in the proof of Theorem
(14), the expansion (5.1), and then equating the coefficients of $z^{6n-2}$.

The asymptotic form of (b) is given by
$$S_{6n}(e^{i\pi/3})\sim {4 \over {(6n-1)!}}\int_0^\infty [e^{-t}\cosh^2(e^{i\pi/3}/2)+(-1)^n
e^{it\exp(i\pi/3)}\cos^2(t/2)]t^{6n-1}dt.$$
The integral may be evaluated in terms of $\Gamma(6n)=(6n-1)!$.  Then after a considerable
number of steps of simplification, we obtain the stated form.  \qed

{\it Remark}.  Based upon numerical investigation, a better approximation to $S_{6n}(e^{i\pi/3})$ 
is given by
$$S_{6n}(e^{i\pi/3})\sim 6+2\left(-{1 \over {18}}\right)^n.$$

\medskip
\centerline{\bf Discussion:  other Bernoulli relations}
\medskip

{\bf Proposition 2}.  For integers $k\geq 0$,
$$B_{2k+2}={{(2k+2)(2k+1)} \over {4(2^{2k+2}-1)}}\int_0^1 E_{2k}(x)dx.$$

{\it Proof}.  We integrate the generating function (1.4), so that
$${z \over 2}\tanh\left({z \over 2}\right)=\sum_{n=0}^\infty {z^n \over {n!}} \int_0^1 E_n(x)dx.$$
We then compare with the generating function
$$\tanh x=\sum_{k=1}^\infty {{2^{2k}(2^{2k}-1)} \over {(2k)!}}B_{2k}x^{2k-1}, ~~~~~~|x|<{\pi \over 2},$$ 
and the result follows.  \qed

Similarly, if we integrate the generating function (2.1) for Bernoulli polynomials, we obtain
the well known relations $\int_0^1 B_0(x)dx=1$ and for $n>0$, $\int_0^1 B_n(x)dx=0$.  

Lehmer obtained the recurrences \cite{carlitz,lehmer}
$$\sum_{k=0}^n{{6n+3} \choose {6k}}B_{6k}=2n+1$$
and
$$\sum_{k=0}^n{{6n+5} \choose {6k+2}}B_{6k+2}={1 \over 3}(6n+5).$$
In terms of the Riemann zeta function at even integers congruent to $0$ and $2$ modulo 6, we thus have
$$\sum_{k=0}^n{{(6n+3)!} \over {(6n-6k+3)!}}{{(-1)^{k+1}} \over {2^{6k-1}\pi^{6k}}}\zeta(6k)=2n+1$$
and
$$\sum_{k=0}^n{{(6n+5)!} \over {(6n-6k+3)!}}{{(-1)^{k+1}} \over {2^{6k+1}\pi^{6k+2}}}\zeta(6k+2)
={1 \over 3}(6n+5).$$
These expressions then provide a starting point for generalization according to the previous
sections of this paper. 
Furthermore, we have determined the following recurrences, which may then be expressed in
terms of sums of values $\zeta(6k+q)$, being special cases of $\zeta(6k+q,a)$, $\Phi(z,6k+q,a)$,
and of other functions.  
{\newline \bf Proposition 3}.  
$$\sum_{k=0}^n{{6n+7} \choose {6k+4}}B_{6k+4}=-\left(n+{7 \over 6}\right),$$
$$\sum_{k=0}^n{{6n+9} \choose {6k+6}}B_{6k+6}=2(n+1),$$
$$\sum_{k=0}^n{{6n+11} \choose {6k+8}}B_{6k+8}=-{1 \over 2}(6n+11)(n+1),$$
$$\sum_{k=0}^n{{6n+13} \choose {6k+10}}B_{6k+10}={1 \over {60}}(6n+13)(18n^2+81n+100),$$
$$\sum_{k=0}^n{{6n+15} \choose {6k+12}}B_{6k+12}=-{1 \over {420}} (648n^4+6156n^3+22266n^2+36765n+24185)(n+2)(n+1),$$
$$\sum_{k=0}^n{{6n+17} \choose {6k+14}}B_{6k+14}={1 \over {8400}}(6n+17)(n+2)(n+1)$$
$$\times(1944n^5+23652n^4+116046n^3+288423n^2+366675n+196000),$$
and
$$\sum_{k=0}^n{{6n+19} \choose {6k+16}}B_{6k+16}=-{1 \over {55440}}(6n+19)(n+2)(n+1)$$
$$\times
(11664n^7+209952n^6+1623240n^5+6998400n^4+182132n^3+28919198n^2+25568993n+10026324).$$
These successive sums may be obtained from the initial ones by shifting the summation index $k$
and upper limit $n$.

Owing to the particular binomial coefficient in the summand in these relations, they are closely
connected with the properties of various $_6F_5$ hypergeometric functions when integral
representations are applied.  Similarly, other $_6F_5$ functions appear if we use a summation
representation such as 
$${B_{2k} \over {2k}}=2\sum_{m=1}^\infty {m^{2k-1} \over {e^{2\pi m}-1}}-{{2\pi} \over k}[1+(-1)^k]
\sum_{m=1}^\infty {{m^{2k}e^{2\pi m}} \over {(e^{2\pi m}-1)^2}}.$$

\bigskip
\centerline{\bf Acknowledgements}

I thank M. C. Lettington for useful discussions and the Cardiff University School of
Mathematics for its hospitality during a Fulbright Commission-sponsored visit during
May-June 2014.

\pagebreak


\begin{thebibliography}{99}
\bibitem{nbs}M. Abramowitz and I. A. Stegun,
{Handbook of Mathematical Functions, Washington, National Bureau of Standards (1964).}
\bibitem{armitage}J. V. Armitage and W. F. Eberlein,
{Elliptic Functions, Cambridge University Press (2006).}
\bibitem{basu}A. Basu and T. M. Apostol,
{A new method for investigating Euler sums, Ramanujan J. {\bf 4}, 397-419 (2000).}
\bibitem{basoco}M. A. Basoco,
{On the trigonometric expansion of elliptic functions, Bull. Amer. Math. Soc. {\bf 37}, 117-124 (1931).}
\bibitem{boutroux}P. Boutroux,
{Recherches sur les transcendantes de M. Painlev\'{e} et l'\'{e}tude asymptotique des \'{e}quations diff\'{e}rentielles du second ordre, Ann. Sci. \'{E}cole Norm. Sup., S\'{e}rie 3 {\bf 30}, 255-375 (1913). (suite) (3) {\bf 31} 99-159 (1914).}
\bibitem{cahill}N. D. Cahill et al.,
{Fibonacci determinants, College Math. J. {\bf 33}, 221-225 (2002).}
\bibitem{carlitz}L. Carlitz,
{Bernoulli numbers, Fib. Q. {\bf 6}, 71-85 (1968).}
\bibitem{carlitz62}L. Carlitz,
{The coefficients of the lemniscate function, Math. Comp. {\bf 16}, 475-478 (1962).}
\bibitem{cofpla2002}M. W. Coffey,
{A set of identities for a theta function at unit argument, Phys. Lett. A {\bf 300}, 367-369 (2002).}
\bibitem{coffeyst1}M. W. Coffey,
{Series representations of the Riemann and Hurwitz zeta functions and series and integral representations of the first Stieltjes constant arXiv:1106.5147 (2011).}
\bibitem{dienstfrey}A. Dienstfrey and J. Huang,
{Integral representations for elliptic functions, J. Math. Anal. Appl. {\bf 316}, 142-160 (2006);
arXiv:0409216v1 (2004).}
\bibitem{frame}J. S. Frame,
{The Hankel power sum matrix inverse and the Bernoulli continued fraction, Math. Comp., {\bf 33}, 
815-826 (1979).} 
\bibitem{grad}I. S. Gradshteyn and I. M. Ryzhik,
{Table of Integrals, Series, and Products, Academic Press, New York (1980).}
\bibitem{grosswald}E. Grosswald,
{Die Werte der Riemannschen Zetafunktion an ungeraden Argumentstellen, Nachr. Akad. Wiss. Gottingen
Math.-Phys. Kl. II, 9-13 (1970).}
\bibitem{gun}S. Gun, M. R. Murty, and P. Rath,
{Transcendental values of certain Eichler integrals, Bull. London Math. Soc. {\bf 43}, 939-952 (2011).}
\bibitem{hurwitz}A. Hurwitz,
{Entwickelungskoeffizienten der lemniscatischen Funktionen, Math. Ann. {\bf 51}, 196-226 (1899).}
\bibitem{ireland}K. Ireland and M. Rosen,
{A Classical Introduction to Modern Number Theory, 2nd ed., Springer (1990).}
\bibitem{lehmer}D. H. Lehmer,
{Lacunary recurrence formulas for the numbers of Bernoulli and Euler, Ann. Math. {\bf 36}, 637-649 (1935).}
\bibitem{mcl1}M. C. Lettington, 
{Fleck's congruence, associated magic squares and a zeta identity, Funct. Approx. Comment. Math. 
{\bf 45}, 165-205 (2011).}
\bibitem{mcl2}M. C. Lettington, 
{A trio of Bernoulli relations, their implications for the Ramanujan polynomials and the
special values of the Riemann zeta function, Acta Arith. {\bf 158}, 1-31 (2013).}
\bibitem{underwood}R. S. Underwood,
{An expression for the summation $\sum_{m=1}^n m^p$, Amer. Math. Monthly {\bf 35}, 424-428 (1928).}
\bibitem{walker}P. Walker,
{On the functional equations satisfied by modular functions, Mathematika {\bf 25}, 185-190 (1978).}
\bibitem{walkerbk}P. Walker,
{Elliptic Functions A Constructive Approach, John Wiley (1996).}
\bibitem{williams}G. T. Williams,
{A new method of evaluating $\zeta(2n)$, Amer. Math. Monthly {\bf 60}, 19-25 (1953).}
\bibitem{zagier}D. Zagier,
{A Kronecker limit formula for real quadratic fields, Math. Ann. {\bf 213}, 153-184 (1975).}
\end{thebibliography}
\end{document}